\documentclass[12pt]{article}
\usepackage{amsmath,amsthm}

\def\be{\begin{equation}}
\def\ee{\end{equation}}
\def\br{\begin{eqnarray*}}
\def\er{\end{eqnarray*}}

\begin{document}
\title{The hypersurfaces with  conformal normal
 Gauss map in $H^{n+1}$ and $S_{1}^{n+1}$
\thanks {Mathematics Subject Classification(2000)  53C42, 53A10}
\thanks { Key words.  the fourth fundamental form, conformal
normal Gauss map, generalized Gauss map, duality property, de Sitter
Gauss map.}
\thanks {This work is supported by the starting foundation of science
research of Shandong University.}}
\author{ Shuguo Shi  }
\date{}
\maketitle
\begin{center}
\begin{minipage}{120mm}

\noindent {\small{\bf Abstract.}
   In this paper  we introduce the fourth fundamental form for the
   hypersurfaces in $H^{n+1}$ and the space-like hypersurfaces in
   $S_{1}^{n+1}$ and discuss the
   conformality of the normal Gauss maps of the hypersurfaces in
   $H^{n+1}$ and $S_{1}^{n+1}$. Particularly,
   we discuss the surfaces with    conformal normal Gauss maps in
    $H^{3}$ and $S_{1}^{3}$ and prove a duality property. We give
    the Weierstrass representation formula for the space-like
    surfaces in $S_{1}^{3}$ with  conformal normal Gauss maps. We
    also state the similar results for the time-like surfaces in
    $S_{1}^{3}.$}
\end{minipage}
\end{center}

\section{Introduction}
It is well known that the classical Gauss map has played an
important role in the study of the surface theory in $R^{3}$ and has
been generalized to the submanifold of arbitrary dimension and
codimension immersed  into the space forms with constant sectional
curvature( see [15]in detail).

Particularly, for the $n$-dimensional submanifold $x:M\rightarrow V$
in space $V$ with constant sectional curvature, Obata[13] introduced
the generalized Gauss map which assigns to each point $p$ of $M$ the
totally geodesic $n$-subspace of $V$ tangent to $x(M)$ at $x(p)$. He
defined the third fundamental form of the  submanifold  in constant
curvature space as the pullback of the metric of the set of all the
totally geodesic $n$-subspaces in $V$ under the generalized Gauss
map. He derived a relationship among the Ricci form of the immersed
submanifold and  the first, the second and the third fundamental
forms of the immersion. Meanwhile, Lawson[10] discussed the
generalized Gauss map of the immersed  surfaces in $S^{3}$ and prove
a duality property between the minimal surfaces in $S^{3}$ and their
generalized Gauss map image.

Epstein[4] and Bryant[3] defined the hyperbolic Gauss map for the
surfaces in $H^{3}$ and Bryant[3] obtained a Weierstrass
representation formula for the constant mean curvature one  surfaces
with conformal hyperbolic Gauss map. Using the Weierstrass
representation formula, Bryant also studied the properties of
constant mean curvature one surfaces. Using the hyperbolic Gauss
map, G{\'{a}}lvez and Mart{{\scriptsize {\'l}}}nez and
Mil{\'{a}}n[6] studied the flat surfaces in $H^{3}$ with  conformal
hyperbolic Gauss map with respect to the second conformal structure
on surfaces (see [7] for the definition) and obtained a Weierstrass
representation formula for such as surfaces.

Kokubu[8] considered the $n$-dimensional hyperbolic space $H^{n}$ as
a Lie group $G$ with a left-invariant metric and defined the normal
Gauss maps of the surfaces which assigns to each  point of the
surface the tangent plane translated to the Lie algebra of $G.$ He
also gave a Weierstrass representation formula for minimal surfaces
in $H^{n}.$ On the other hand, G{\'{a}}lvez and Mart{{\scriptsize
{\'l}}}nez[5] studied the properties of the Gauss map of a surface
$\Sigma$ immersed into the Euclidean 3-space $R^{3}$ by using the
second  conformal structure on surface and obtained the Weierstrass
representation formula for the surfaces with prescribed Gauss map.
Motivated by their work, the author[16] gave a Weierstrass
representation formula for the surfaces with prescribed normal Gauss
map and Gauss curvature in $H^{3}$ by using the second conformal
structure on surfaces. From this, the surfaces whose normal Gauss
maps are conformal have been found and the translational surfaces
with conformal normal Gauss maps locally are given. In [17], the
author classified locally the ruled surfaces with conformal normal
Gauss maps within the Euclidean ruled surfaces and studied some
global properties of the ruled surfaces and translational surfaces
with conformal normal Gauss maps.

Aiyama and Akutagawa [1] defined the normal Gauss map for the
space-like surfaces in the de Sitter 3-space $S_{1}^{3}$ and gave
the Weierstrass representation formula for the space-like surfaces
in $S_{1}^{3}$ with prescribed mean curvature and normal Gauss map.

The purpose of this paper is to study the conformality of the normal
Gauss maps for the hypersurfaces in $H^{n+1}$ and the space-like
hypersurfaces in $S_{1}^{n+1}$ and to prove a duality property
between the surfaces in $H^{3}$ and the space-like surfaces in
$S_{1}^{3}$ with  conformal normal Gauss maps. The rest of this
paper is organized as follows. In the second section, we describe
the generalized definition of the normal Gauss map for the
hypersurfaces in $H^{n+1}$ and the space-like hypersurfaces in
$S_{1}^{n+1}$ (cf.[1][8]). The third section introduces the fourth
fundamental form for the hypersurfaces in $H^{n+1}$ and
$S_{1}^{n+1}$  and obtains a relation among the first, the second,
the third and the fourth fundamental forms of the hypersurfaces. As
a application, we discuss the conformality of the normal Gauss map
for the hypersurfaces in $H^{n+1}$ and the space-like hypersurfaces
in $S_{1}^{n+1}$. By means of the generalized Gauss map of the
surfaces in $H^{3}$ and $S_{1}^{3},$ the fourth one proves a duality
property between the surfaces in $H^{3}$ and the space-like surfaces
in $S_{1}^{3}$ with
 conformal normal Gauss maps.  The fifth one gives the Weierstrass
representation formula for the space-like surfaces in $S_{1}^{3}$
with  conformal normal Gauss map and the sixth one
 derives the PDE for the
space-like graphs in $S_{1}^{3}$ with  conformal normal Gauss map
and classifies locally the translational surfaces and the Euclidean
ruled surfaces in  $S_{1}^{3}$ with conformal normal Gauss map. In
the last section, we state the similar results for time-like
surfaces in $S_{1}^{3}$ with conformal normal Gauss map.\vspace{3mm}

{\bf Acknowledgement}\hspace{1mm} The author would like to express
his sincere gratitude to Prof. Detang Zhou for his enthusiastic
encouragement, support and valuable help as well as for his
significant suggestions and heuristic discussions with the author
and for his providing  the author with Omori and Yau's paper
[14][18].

\section{Preliminaries}\setcounter{equation}{0}
Take the upper half-space models of the hyperbolic space
$H^{n+1}(-1)$ and  the de Sitter space $S_{1}^{n+1}(1)$
$$ R_{+}^{n+1}=\left \{(x_{1},x_{2},\cdots, x_{n+1})\in
R^{n+1}|x_{n+1}>0\right\}$$ with respectively the Riemannian
metric $ds^{2}={\frac
{1}{x_{n+1}^{2}}}(dx_{1}^{2}+dx_{2}^{2}+\cdots +dx_{n+1}^{2})$
 and the Lorentz metric $ds^{2}={\frac
{1}{x_{n+1}^{2}}}(dx_{1}^{2}+dx_{2}^{2}+\cdots
+dx_{n}^{2}-dx_{n+1}^{2})$ (cf.[1]).

Let $M$ be a $n$-dimensional  Riemannian manifold and
$x:M^{n}\rightarrow H^{n+1}$($resp.$ $x:M^{n}\rightarrow
S_{1}^{n+1}$) be an immersed hypersurface ($resp.$ space-like
hypersurface) with the local coordinates $u_{1},u_{2},\cdots,u_{n}.$
 In this paper,we agree with the following ranges of indices: $1\leq
i,j,k,\cdots\leq n$ and $1\leq A,B,C,\cdots\leq n+1.$ The first and
the second fundamental forms are given, respectively, by $
\textrm{I}=g_{ij}du_{i}du_{j}$ and $\textrm{II}=h_{ij}du_{i}du_{j}.$
 The unit normal vector ($resp.$ time-like unit normal vector) of $x(M)$ is
$N=x_{n+1}\eta_{1}{\frac {\partial}{\partial x_{1}}}+
x_{n+1}\eta_{2}{\frac {\partial}{\partial x_{2}}}+\cdots+
x_{n+1}\eta_{n+1}{\frac {\partial}{\partial x_{n+1}}},$ where
$\eta_{1}^{2}+\eta_{2}^{2}+\cdots +\eta_{n+1}^{2}=1$($resp.$ $
\eta_{1}^{2}+\eta_{2}^{2}+\cdots +\eta_{n}^{2}-\eta_{n+1}^{2}=-1$).

We have the Weingarten formula
$$ {\frac {\partial \eta_{A}}{\partial u_{k}}}={\frac
{1}{x_{n+1}}}\left (\eta_{n+1}{\frac {\partial x_{A}}{\partial
u_{k}}}-g^{jl}h_{kl}{\frac {\partial x_{A}}{\partial u_{j}}}\right
)$$$$\left(resp. {\frac {\partial \eta_{A}}{\partial
u_{k}}}={\frac {1}{x_{n+1}}}\left (\eta_{n+1}{\frac {\partial
x_{A}}{\partial u_{k}}}+g^{jl}h_{kl}{\frac {\partial
x_{A}}{\partial u_{j}}}\right )\right).$$

Identitying $H^{n+1}$ and $S_{1}^{n+1}$ with the Lie group (cf.[8])
\[
 {G}=\left \{\left (
\begin{array}{lcccr}
1 & 0 & \cdots & 0 & \log x_{n+1}\\
0 & x_{n+1} & \cdots & 0 & x_{1}\\
\vdots & \vdots & \ddots & \vdots & \vdots\\
0 & 0 & \cdots & x_{n+1} & x_{n}\\
0 & 0 & \cdots & 0 & 1
\end{array} \right):(x_{1},x_{2},\cdots,x_{n+1})\in R_{+}^{n+1} \right
\},
\]
the multiplication is defined as the matrix multiplication and the
identity is $e=(0,0,\cdots,0,1).$ The Riemannian metric of $H^{n+1}
$ and the Lorentz metric of $S_{1}^{n+1}$ are left-invariant and
$\widetilde {X}_{1}=x_{n+1}{\frac {\partial }{\partial
x_{1}}},\widetilde {X}_{2}=x_{n+1}{\frac {\partial }{\partial
x_{2}}},\cdots,\widetilde {X}_{n+1}=x_{n+1}{\frac {\partial
}{\partial x_{n+1}}}$ are the left-invariant unit orthonormal vector
fields. Now, the unit normal vector ($resp.$ time-like unit normal
vector) field of $x(M)$ can be written as $N=\eta_{1}{\widetilde
X_{1}}+\eta_{2}{\widetilde X_{2}}+\cdots+\eta_{n+1}{\widetilde
X_{n+1}}.$ Left translating $N$ to $T_{e}(R_{+}^{n+1}),$ we obtain
$$ \widetilde N:M\rightarrow S^{n}(1)\subset
T_{e}(R_{+}^{n+1})(resp.\widetilde N:M\rightarrow H^{n}(-1)\subset
T_{e}(R_{+}^{n+1})),$$
$$\widetilde N=L_{x^{-1}*}(N)=\eta_{1}{\frac {\partial}{\partial x_{1}}}(e)+
\eta_{2}{\frac {\partial}{\partial x_{2}}}(e)+\cdots+
\eta_{n+1}{\frac {\partial}{\partial x_{n+1}}}(e).$$ Call
$\widetilde N$ the normal Gauss map of the immersed hypersurface
$x:M\rightarrow H^{n+1}$($resp.$ space-like hypersurface
$x:M\rightarrow S_{1}^{n+1}$)(cf.[1][8]).

\section{The fourth fundamental form}\setcounter{equation}{0}
{\bf Definition.} {\sl Let $M$ be a $n$-dimensional Riemannian
manifold. Call $\textrm{IV}=\langle {d\widetilde N},{d\widetilde
N}\rangle$ the fourth fundamental form of the immersed hypersurface
$x:M\rightarrow H^{n+1}$($resp.$ space-like hypersurface
$x:M\rightarrow S_{1}^{n+1}$), where the scalar product
$\langle\cdot,\cdot\rangle$ is induced by the Euclidean metric of
$R^{n+1}$ ($resp.$ the Lorentz-Minkowski metric of
$L^{n+1}$).}\vspace{3mm}

\noindent{\bf THEOREM 3.1.} {\sl Let $M$ be a $n$-dimensional
Riemannian manifold with  Ricci form $Ric.$ Let $x:M\rightarrow
H^{n+1}$($resp.$ $x:M\rightarrow S_{1}^{n+1}$) be an immersed
hypersurface ($resp.$space-like hypersurface) with  mean curvature
$H={\frac {1}{n}}tr(\textrm{II}).$ Then
\begin{equation}
\textrm{IV}=\eta_{n+1}^{2}\textrm{I}-2\eta_{n+1}\textrm{II}+\textrm{III}
\end{equation}
\begin{equation}
(resp.\textrm{IV}=\eta_{n+1}^{2}\textrm{I}+2\eta_{n+1}\textrm{II}+\textrm{III}),
\end{equation}
 where
$\textrm{III}=nH\textrm{II}-(n-1)\textrm{I}-Ric$  $
(resp.\textrm{III}=nH\textrm{II}-(n-1)\textrm{I}+Ric)$ is Obata's
third fundamental form of $x(M)$ }(see [13]). \vspace{3mm}

{\sl Proof.} At first we prove the Theorem for $H^{n+1}$. Choose the
normal coordinates $u_{1},u_{2},\cdots,u_{n}$  near $p\in M.$ By the
Weingarten formula, we get
$$
\textrm{IV}=\langle{d\widetilde N},{d\widetilde N}\rangle= {\frac
{\partial \eta_{A}}{\partial u_{i}}}{\frac {\partial
\eta_{A}}{\partial u_{j}}}du_{i}du_{j} $$$$={\frac
{1}{x_{n+1}^{2}}}\left (\eta_{n+1}{\frac {\partial x_{A}}{\partial
u_{i}}}-h_{ik}{\frac {\partial x_{A}}{\partial u_{k}}}\right
)\left (\eta_{n+1}{\frac {\partial x_{A}}{\partial
u_{j}}}-h_{jl}{\frac {\partial x_{A}}{\partial u_{l}}}\right
)du_{i}du_{j}$$
\begin{equation}
=(\eta_{n+1}^{2}\delta_{ij}-2\eta_{n+1}h_{ij}+h_{ik}h_{jk})du_{i}du_{j}.
\end{equation}
$\textrm{III}=h_{ik}h_{jk}du_{i}du_{j}$ is the third fundamental
form [13] and by the Gauss equation,
$\textrm{III}=nH\textrm{II}-(n-1)\textrm{I}-Ric.$ (3.1) is proved.

Next, similar to the above proof,  for $S_{1}^{n+1},$ we have
\begin{equation}
\textrm{IV}=(\eta_{n+1}^{2}\delta_{ij}+2\eta_{n+1}h_{ij}+h_{ik}h_{jk})du_{i}du_{j}.
\end{equation}
Similar to the proof of (3.1), we can prove (3.2). \vspace{3mm}

Next, we consider the applications of these formulas (3.1)$-$(3.4).
In the following of this paper, that the normal Gauss map is
conformal means that the fourth fundamental form is proportional to
the second fundamental form, i.e. $\textrm{IV}=\rho \textrm{II}$ for
some smooth function $\rho$ on $M.$ \vspace{3mm}

 \noindent{\bf THEOREM 3.2.} {\sl Let $M$ be a $n$-dimensional Riemannian manifold and
 $x:M\rightarrow H^{n+1}$($resp.$ $x:M\rightarrow S_{1}^{n+1}$) be an
 immersed hypersurface ($resp.$
space-like hypersurface) without  umbilics. Then the normal Gauss
map of $x(M)$ is conformal  if and only if at each  point of $M$,
there exists exactly two distinct principal curvatures  and the
sectional curvature $R(X\wedge Y)=-1+\eta_{n+1}^{2}$($resp.R(X\wedge
Y)=1-\eta_{n+1}^{2}$), where the  vectors $X$ and $Y$ belong to
different principal direction spaces.}\vspace{3mm}

{\sl Proof.} The case of $H^{n+1}$. For any point $p\in M$, let
$\left\{e_{1},e_{2},\cdots,e_{n}\right \}$ be a local frame field so
that $(h_{ij})$ is diagonalized at this point, i.e.
$h_{ij}(p)=\lambda_{i}\delta_{ij}.$ By $ \textrm{IV}=\rho
\textrm{II}$ and (3.3),  we get, for $i=1,2,\cdots,n,$ that
\begin{equation}
\eta_{n+1}^{2}-2\eta_{n+1}\lambda_{i}+\lambda_{i}^{2}=\rho\lambda_{i},
\end{equation}
i.e.
\begin{equation}
\lambda_{i}^{2}-(\rho+2\eta_{n+1})\lambda_{i}+\eta_{n+1}^{2}=0.
\end{equation}
Because $x(M)$ has no umbilics, the equation (3.6) with respect to
$\lambda_{i}$ has exactly two distinct solutions $\lambda$ and $\mu$
and $\lambda\mu=\eta_{n+1}^{2}$. By the Gauss equation, one may
prove $R(X\wedge Y) =-1+\lambda\mu=-1+\eta_{n+1}^{2}.$

Conversely, choose  the local tangent frame $\left \{
e_{1},e_{2},\cdots, e_{n}\right\}$ and the dual frame
$\{\omega_{1},\omega_{2}, \cdots,\omega_{n}\}$ near $p$, such that
$h_{ij}=0, i\not= j$ and $ h_{11}=h_{22}=\cdots=h_{rr}=\lambda\not=
\mu=h_{r+1r+1}=\cdots=h_{nn}.$ Then $\eta_{n+1}^{2}=\lambda\mu.$ By
(3.3),
\begin{eqnarray*}
\textrm{IV} &=&
(\eta_{n+1}^{2}-2\eta_{n+1}\lambda+\lambda^{2})(\omega_{1}^{2}+\cdots+\omega_{r}^{2})\\
&&+(\eta_{n+1}^{2}-2\eta_{n+1}\mu+\mu^{2})(\omega_{r+1}^{2}+\cdots+\omega_{n}^{2})\\
&=&(\mu-2\eta_{n+1}+\lambda)\lambda(\omega_{1}^{2}+\cdots+\omega_{r}^{2})\\
&&+(\lambda-2\eta_{n+1}+\mu)\mu(\omega_{r+1}^{2}+\cdots+\omega_{n}^{2})\\
&=&(\lambda-2\eta_{n+1}+\mu)\textrm{II}.
\end{eqnarray*}
The sufficiency has been proved for $H^{n+1}$. Similarly, we can
prove Theorem 3.2 for $S_{1}^{n+1}.$\vspace{3mm}

{\sl Remark.} By (3.5), we know that the normal Gauss maps of all
totally umbilics hypersurfaces except the totally geodesic
hyperspheres in $H^{n+1}$  are conformal. Similarly, for the
space-like hypersurfaces in $S_{1}^{n+1}$, since $\eta_{n+1}\not
=0,$ the normal Gauss maps of all totally umbilic space-like
hypersurfaces except totally geodesic space-like hypersurfaces are
conformal.\vspace{3mm}

For $H^{3}$ and $S_{1}^{3}$, by Theorem 3.2, we immediately
get\vspace{3mm}

\noindent{\bf THEOREM 3.3.} {\sl Let $M$ be a 2-dimensional
Riemannian manifold and $x:M\rightarrow H^{3}$($resp.$
$x:M\rightarrow S_{1}^{3}$) be an immersed surface ($resp.$
space-like surface) without  umbilics. Then the normal Gauss map of
$x(M)$ is conformal if and only if the Gauss curvature
$K=-1+\eta_{3}^{2}$($resp.$$K=1-\eta_{3}^{2}$).}\vspace{3mm}

{\sl Remark.} In [16][17], we assume that the second fundamental
form is positive definite and induces the conformal structure on the
surfaces in  $H^{3}$. Here, the assumption with respect to the
positive definite second fundamental form is dropped.\vspace{3mm}

\noindent{\bf THEOREM 3.4.} {\sl Let $M$ be a $n$-dimensional
Einstein manifold and $x:M\rightarrow H^{n+1}$($resp.$
$x:M\rightarrow S_{1}^{n+1}$) be an immersed hypersurface ($resp.$
space-like hypersurface)  with the non-degenerate second fundamental
form and without umbilics. If the normal Gauss map of $x(M)$ is
conformal map,i.e. $\textrm{IV}=\rho \textrm{II},$ then $n=2$ and
$\rho=2(H-\eta_{3})$ ($resp. \rho=2(H+\eta_{3})$).}\vspace{3mm}

{\sl Proof.} We only prove the Theorem for $H^{n+1}$. $ M $ is an
Einstein manifold, so $Ric={\frac {S}{n}}I,$ where $S$ is the scalar
curvature of $M$. (3.1) becomes
$$ \left (\eta_{n+1}^{2}-(n-1)-{\frac {S}{n}}\right )\textrm{I}+
(nH-2\eta_{n+1}-\rho)\textrm{II}=0.$$ Because $x(M)$ has no
umbilics, we have
$$nH=2\eta_{n+1}+\rho.$$
By Theorem 3.2 and its proof, we assume that
$\lambda_{1}=\cdots=\lambda_{r}=\lambda \not=\mu=
\lambda_{r+1}=\cdots\lambda_{n},$ then
$$r\lambda+(n-r)\mu=2\eta_{n+1}+\rho .$$
By (3.6),
$$ \lambda+\mu=2\eta_{n+1}+\rho.$$
So $(r-1)\lambda+(n-r-1)\mu=0.$ By Theorem 3.2, $\lambda$ and $\mu$
have same signature. So $r=1$ and $n=2$.  Hence $\rho=2H-2\eta_{3}.$

\section{A duality for the  surfaces in $H^{3}$ and
$S_{1}^{3}$ with  conformal normal Gauss maps}
\setcounter{equation}{0} Let $L^{4}$ be the  Minkowski 4-space
 with the canonical coordinates $X_{0},X_{1}, X_{2},X_{3}$ and
the Lorentz-Minkowski scalar product
$-X_{0}^{2}+X_{1}^{2}+X_{2}^{2}+X_{3}^{2}.$  The  Minkowski model of
$H^{3}$ is given by
$$H^{3}=\{(X_{0},X_{1},X_{2},X_{3})|-X_{0}^{2}+
X_{1}^{2}+X_{2}^{2}+X_{3}^{2}=-1,X_{0}>0\}$$ and is identified
with the upper half-space model $R_{+}^{3}$ of $H^{3}$ by
$$ (x_{1}, x_{2},x_{3})=\left ({\frac {X_{1}}{X_{0}-X_{3}}},
{\frac {X_{2}}{X_{0}-X_{3}}},{\frac {1}{X_{0}-X_{3}}}\right).$$
Accordingly, the space-like normal vector of the surface in the
Minkowski model of $H^{3}$ is $N=N_{0}{\frac {\partial}{\partial
X_{0}}}+ N_{1}{\frac {\partial}{\partial X_{1}}}+N_{2}{\frac
{\partial}{\partial X_{2}}}+N_{3}{\frac {\partial}{\partial
X_{3}}},$ where
\begin{eqnarray*}
N_{0}&=&{\frac {X_{1}}{X_{0}-X_{3}}}\eta_{1}+{\frac
{X_{2}}{X_{0}-X_{3}}}\eta_{2}+{\frac
{1-X_{0}(X_{0}-X_{3})}{X_{0}-X_{3}}}\eta_{3},\\
N_{1}&=& \eta_{1}-X_{1}\eta_{3}, N_{2}=\eta_{2}-X_{2}\eta_{3},\\
N_{3}&=&{\frac {X_{1}}{X_{0}-X_{3}}}\eta_{1}+{\frac
{X_{2}}{X_{0}-X_{3}}}\eta_{2}+{\frac
{1-X_{3}(X_{0}-X_{3})}{X_{0}-X_{3}}}\eta_{3}.
\end{eqnarray*}
We get
\begin{equation}
\eta_{3}={\frac {N_{0}-N_{3}}{X_{3}-X_{0}}}.
\end{equation}

The Minkowski model of the de Sitter 3-space is defined as
$$ S_{1}^{3}=\{(X_{0},X_{1},X_{2},X_{3})|-X_{0}^{2}+
X_{1}^{2}+X_{2}^{2}+X_{3}^{2}=1\}\simeq S^{2}\times R$$ and can be
divided into three components as follows(cf. [1]),
$$S_{-}=\{X\in S_{1}^{3}|X_{0}-X_{3}<0\}\simeq R^{3},$$
$$S_{0}=\{X\in S_{1}^{3}|X_{0}-X_{3}=0\}\simeq S^{1}\times R,$$
$$S_{+}=\{X\in S_{1}^{3}|X_{0}-X_{3}>0\}\simeq R^{3}.$$
Identify $S_{-}$ and $S_{+}$ with the  upper half-space model
$R_{+}^{3}$ of the de Sitter 3-space by (cf. [1])
$$ (x_{1}, x_{2},x_{3})=\left ({\frac {X_{1}}{|X_{0}-X_{3}|}},
{\frac {X_{2}}{|X_{0}-X_{3}|}},{\frac
{1}{|X_{0}-X_{3}|}}\right).$$ For the space-like surface
$X:M\rightarrow S_{1}^{3}$, let $U_{-}=X^{-1}(S_{-})$ and
$U_{+}=X^{-1}(S_{+}),$ then $U_{-}\cup U_{+}$ is the open dense
subset of $M.$ On $U_{-}\cup U_{+}$, the time-like unit normal
vector is $N=N_{0}{\frac {\partial}{\partial X_{0}}}+ N_{1}{\frac
{\partial}{\partial X_{1}}}+N_{2}{\frac {\partial}{\partial
X_{2}}}+N_{3}{\frac {\partial}{\partial X_{3}}},$ where
\begin{eqnarray*}
N_{0}&=&{\frac {X_{1}}{X_{0}-X_{3}}}\eta_{1}+{\frac
{X_{2}}{X_{0}-X_{3}}}\eta_{2}-{\frac
{1+X_{0}(X_{0}-X_{3})}{X_{0}-X_{3}}}\eta_{3},\\
N_{1}&=& \eta_{1}-X_{1}\eta_{3}, N_{2}=\eta_{2}-X_{2}\eta_{3},\\
N_{3}&=&{\frac {X_{1}}{X_{0}-X_{3}}}\eta_{1}+{\frac
{X_{2}}{X_{0}-X_{3}}}\eta_{2}-{\frac
{1+X_{3}(X_{0}-X_{3})}{X_{0}-X_{3}}}\eta_{3}.
\end{eqnarray*}
We get
\begin{equation}
\eta_{3}={\frac {N_{0}-N_{3}}{X_{3}-X_{0}}}.
\end{equation}

{\sl Remark.} In [1], the normal Gauss map of the space-like surface
$X:M\rightarrow S_{1}^{3}$ is defined globally on $M.$ Because of
the density of $U_{-}$ and $U_{+}$ in $M,$ in this paper, we may
consider that the normal Gauss map of the space-like surface
$X:M\rightarrow S_{1}^{3}$ is defined on $U_{-}$ and
$U_{+}.$\vspace{3mm}

Let $X:M\rightarrow H^{3}$($resp.$ $X:M\rightarrow S_{1}^{3}$) be
an  immersed surface($resp.$ space-like surface). Parallel
translating the space-like ($resp.$ time-like) unit normal vector
$N$ to the origin of $L^{4}$, one gets the map $N: M\rightarrow
S_{1}^{3}$ ($resp.N: M\rightarrow H^{3}$) which is usually called
generalized Gauss map of $X:M\rightarrow H^{3}$($resp.$
$X:M\rightarrow S_{1}^{3}$). The generalized Gauss map image can
be considered as the surface in  $S_{1}^{3}$($resp.
H^{3}$).\vspace{3mm}

\noindent{\bf THEOREM 4.1}(cf[9]. Prop 3.5). {\sl(1) Let
$X:M\rightarrow H^{3}$ be a 2-dimensional immersed surface. Then its
generalized Gauss map $N:M\rightarrow S_{1}^{3}$ is a branched
space-like immersion into $S_{1}^{3}$ with branch points where
$K=-1.$ And, when $K\not=-1,$ the curvature of $N:M\rightarrow
S_{1}^{3}$ is $K^{*}={\frac {K}{K+1}}$ and the volume element is
$dV_{N}=|K+1|dV_{X}.$

(2) Let $X:M\rightarrow S_{1}^{3}$ be a 2-dimensional space-like
immersed surface. Then its generalized Gauss map $N:M\rightarrow
H^{3}$ is a branched  immersion into $H^{3}$ with branch points
where $K=1.$ And, when $K\not=1,$ the curvature of $N:M\rightarrow
H^{3}$ is $K^{*}={\frac {K}{1-K}}$ and the volume element is
$dV_{N}=|1-K|dV_{X}.$}\vspace{3mm}

{\sl Proof.} In the context of this paper, we  prove (2). For any
$p\in M$, let$\{e_{0},e_{1},e_{2},e_{3}\}$ be the orthonormal frame
near $p$, such that $ e_{0}=X, e_{3}=N.$ Let
$\{\omega_{0},\omega_{1},\omega_{2},\omega_{3}\}$ be the dual frame.
The connection 1-forms  is $\omega_{\alpha}^{\beta},$
$\alpha,\beta=0,1,2,3.$ The coefficients of the second fundamental
form of $X:M\rightarrow S_{1}^{3}$ is given by
$\omega_{i}^{3}=h_{ij}\omega_{j}, h_{ij}=h_{ji}, i,j=1,2.$ The
induced metric of $N :M\rightarrow H^{3}$ is $ds_{*}^{2}=\langle
dN,dN\rangle =h_{ik}h_{jk}\omega_{i}\omega_{j}.$ Choose the local
tangent frame $\{e_{1},e_{2}\}$ near  $p$ , such that
$h_{ij}=\lambda_{i}\delta_{ij}.$ Then
$ds_{*}^{2}=\lambda_{1}^{2}\omega_{1}^{2}+
\lambda_{2}^{2}\omega_{2}^{2}.$ So, when $\lambda_{1}\lambda_{2}\not
= 0,$ i.e. $K\not=1,$ $ N(M)$ is an immersed  surface into $H^{3}$.
Its space-like unit normal vector is $X$ and the second fundamental
form is $\textrm{II}=-\langle
dX,dN\rangle=-\lambda_{1}\omega_{1}^{2}-\lambda_{2}\omega_{2}^{2}.$
By the Gauss equation,  $K^{*}=-1+{\frac
{1}{\lambda_{1}\lambda_{2}}}={\frac {K}{1-K}}.$\vspace{3mm}

By Theorem 3.3, (4.1),(4.2) and Theorem 4.1, we get the following
duality. \vspace{3mm}

\noindent{\bf THEOREM 4.2.} {\sl Let $M$ be a connected
2-dimensional manifold. Let $X:M\rightarrow H^{3}$ be an immersed
surface without umbilics and $K\not= -1$ and let  $N:M\rightarrow
S_{1}^{3}$ be a space-like surface without umbilics and $K\not=1.$
Suppose that $N:M\rightarrow S_{1}^{3}$ is the generalized Gauss map
of $X:M\rightarrow H^{3}$ and vice versa. Then, the normal Gauss map
of $X:M\rightarrow H^{3}$ is conformal if and only if one of
$N:M\rightarrow S_{1}^{3}$ is conformal. And, at this time,
$dV_{N}=\left ({\frac {N_{0}-N_{3}}{X_{3}-X_{0}}}\right
)^{2}dV_{X}.$}\vspace{3mm}

{\sl Remark.} Like [10] for minimal surfaces in $S^{3}$, we call the
generalized Gauss map  $N: M\rightarrow S_{1}^{3}$ the   polar
variety of the immersed surface $ X:M\rightarrow H^{3}$  with
conformal normal Gauss map and vice versa.

\section{Weierstrass representation
formula}\setcounter{equation}{0} In this section, we give the
Weierstrass representation formula for the space-like surfaces in
$S_{1}^{3}$ with  conformal normal Gauss maps. At first, we describe
the normal Gauss map and the de Sitter Gauss map of the space-like
surfaces in $S_{1}^{3}.$ Take the upper half-space model $R_{+}^{3}$
of $S_{1}^{3}.$

The normal Gauss map of the space-like surface $x:M\rightarrow
S_{1}^{3}$ is given by $\widetilde N=\eta_{1}{\frac
{\partial}{\partial x_{1}}}(e)+ \eta_{2}{\frac {\partial}{\partial
x_{2}}}(e)+ \eta_{3}{\frac {\partial}{\partial
x_{3}}}(e):M\rightarrow H^{2}(-1)\subset L^{3}.$  By means of the
stereographic projection from the north pole $(0,0,1)$ of
$H^{2}(-1)$ to the $(x_{1},x_{2})-$plane identified with $C,$ we get
$$g^{S}={\frac {\eta_{1}+i\eta_{2}}{1-\eta_{3}}}:M\rightarrow C\cup
\{\infty\}\backslash \{|z|=1\},$$ which is also called the normal
Gauss map of the space-like surface $x:M\rightarrow S_{1}^{3}.$
$\widetilde N$ can be written as $$\widetilde N=\left (-{\frac
{g+{\bar g}}{|g|^{2}-1}},{\frac {i(g-{\bar g})}{|g|^{2}-1}},{\frac
{1+|g|^{2}}{|g|^{2}-1}}\right ).$$

Next, we describe the definition of the de Sitter  Gauss map for the
space-like surfaces in $S_{1}^{3}$(in [11], it is still called
hyperbolic Gauss map), which is the analogue of Epstein and Bryant's
hyperbolic Gauss map for the surfaces in $H^{3}$(cf [3][4][16]). The
time-like geodesic is either the Euclidean equilateral
half-hyperbola consisting of  two branches which is orthonormal to
the coordinate plane $\{(x_{1},x_{2},0)|(x_{1},x_{2})\in R^{2}\}$ or
the Euclidean straight line which is orthonormal to the above
coordinate plane. For the space-like surface
$x=(x_{1},x_{2},x_{3}):M\rightarrow S_{1}^{3},$ at each point $x\in
M,$ the oriented time-like geodesic in $S_{1}^{3}$ passing through
$x$ with the time-like tangent vector $N$ meets
$\{(x_{1},x_{2},0)|(x_{1},x_{2})\in R^{2}\}\cup\{\infty\}$ two
points. Since the geodesic is oriented, we may speak of one of the
two points as the initial point and the other one as the final
point. Call the final point the image of the de Sitter Gauss map for
$x(M)$ at the point $x.$ Denote the de Sitter Gauss map by $G^{S}.$
On the coordinate plane $\{(x_{1},x_{2},0)|(x_{1},x_{2})\in
R^{2}\}$, we introduce the natural complex coordinate
$z=x_{1}+ix_{2}.$ Using the Euclidean geometry, as similar as done
in the Theorem 5.1 of [16], we get
\begin{equation}
G^{S}=x_{1}+ix_{2}+x_{3}g^{S}.
\end{equation}

Let $x=(x_{1},x_{2},x_{3}):M\rightarrow H^{3}$ be an immersed
surface with  unit normal vector $N=x_{3}\eta_{1}{\frac
{\partial}{\partial x_{1}}}+ x_{3}\eta_{2}{\frac {\partial}{\partial
x_{2}}}+x_{3}\eta_{3}{\frac {\partial}{\partial x_{3}}}.$ By the
duality given in section 4, the generalized Gauss map of
$x:M\rightarrow H^{3}$ is given, when $\eta_{3}>0$, by
\begin{equation}
N=\left ( {\frac {\eta_{1}}{\eta_{3}}}x_{3}-x_{1}, {\frac
{\eta_{2}}{\eta_{3}}}x_{3}-x_{2},{\frac {x_{3}}{\eta_{3}}}\right): M
\rightarrow S_{1}^{3},
\end{equation}
and when $\eta_{3}<0$, by
\begin{equation}
N=\left ( x_{1}-{\frac {\eta_{1}}{\eta_{3}}}x_{3}, x_{2}-{\frac
{\eta_{2}}{\eta_{3}}}x_{3},-{\frac {x_{3}}{\eta_{3}}}\right): M
\rightarrow S_{1}^{3}
\end{equation}
and in the Minkowski model of the de Sitter 3-space, their time-like
unit normal vector is $X: M\rightarrow H^{3}.$ Again by the duality
given in section 4, a straightforward computation shows us that the
normal Gauss map of $N: M\rightarrow S_{1}^{3}$ is given by
 $${\widetilde N}={\frac {\eta_{1}}{\eta_{3}}}{\frac
{\partial}{\partial x_{1}}}(e)+{\frac {\eta_{2}}{\eta_{3}}}{\frac
{\partial}{\partial x_{2}}}(e)+{\frac {1}{\eta_{3}}}{\frac
{\partial}{\partial x_{3}}}(e): M\rightarrow H^{2}(-1).$$ So,
\begin{equation}
g^{S}={\frac {{\frac {\eta_{1}}{\eta_{3}}}+i{\frac
{\eta_{2}}{\eta_{3}}}}{1-{\frac {1}{\eta_{3}}}}}={\frac
{\eta_{1}+i\eta_{2}}{\eta_{3}-1}}=-g^{H},
\end{equation}
where $g^{H}:M\rightarrow C\cup\{\infty\}$ is exactly the normal
Gauss map of $x:M\rightarrow H^{3}$(cf[8][16][17]). From this, we
also prove the Theorem 4.2.

By (5.1)-(5.4) and the Theorem 5.1 of [16], we get that when
$\eta_{3}>0,$ i.e. $|g^{S}|>1,$
\begin{equation}
G^{S}=-G^{H},
\end{equation}
and when $\eta_{3}<0,$ i.e. $|g^{S}|<1,$
\begin{equation}
G^{S}=G^{H},
\end{equation}
where $G^{H}$ is exactly the hyperbolic Gauss map of $x:M\rightarrow
H^{3}$(cf[3][4][16]).

In the following, we write respectively $g^{S}$ and $G^{S}$ as $g$
and $G.$

By (5.2)-(5.6) and the Weierstrass representation for the surfaces
in $H^{3}$ with  conformal normal Gauss map[16], we get the
Weierstrass representation formula for the space-like surfaces in
$S_{1}^{3}$ with  conformal normal Gauss map.\vspace{3mm}

\noindent{\bf THEOREM 5.1.} {\sl Let $M$ be a simply connected
Riemannian surface. Given the map $G:M\rightarrow C\cup\{\infty\}$
and the nonconstant conformal map $g:M\rightarrow
C\cup\{\infty\}\backslash \{|z|=1\}.$

(1) When the holomorphic map $g:M\rightarrow
C\cup\{\infty\}\backslash \{|z|=1\}$ satisfies $|g|>1$ and
\begin{eqnarray}
&&{\frac {G_{z}}{g_{z}}}>0,\\
&&|g|^{2}|G_{\bar z}|>|G_{z}|,\\
&&G_{z{\bar z}}+{\frac {{\bar g}_{\bar z}}{(|g|^{4}-1){\bar
g}}}G_{z}-{\frac {|g|^{2}{\bar g}g_{ z}}{|g|^{4}-1}}G_{\bar z}=0,
\end{eqnarray}
put
\begin{eqnarray}
x_{1}&=& Re\left \{G-{\frac {1+|g|^{2}}{{\bar
g}g_{z}}}G_{z}\right\},\\
x_{2}&=& Im\left \{G-{\frac {1+|g|^{2}}{{\bar
g}g_{z}}}G_{z}\right\},\\
x_{3} &=& {\frac {1+|g|^{2}}{|g|^{2}g_{z}}}G_{z}.
\end{eqnarray}
Then $x=(x_{1},x_{2},x_{3}):M\rightarrow S_{1}^{3}$ is a space-like
surface with de Sitter Gauss map $G$ and  holomorphic normal Gauss
map $g$ and Gauss curvature $K$ satisfying ${\sqrt {1-K}}={\frac
{1+|g|^{2}}{|g|^{2}-1}}.$ And the conformal structure on $M$ is
induced by the negative definite second fundamental form.
Conversely, any surface $x:M\rightarrow S_{1}^{3}$ with ${\sqrt
{1-K}}={\frac {1+|g|^{2}}{|g|^{2}-1}}(=\eta_{3})$ can be given by
(5.10)(5.11)(5.12) and the de Sitter  Gauss map $G$ and the normal
Gauss map $g$ must satisfy (5.7)(5.8)(5.9), where the conformal
structure on $M$ is induced by the negative definite second
fundamental form.}

{\sl(2) When the antiholomorphic map $g:M\rightarrow
C\cup\{\infty\}\backslash \{|z|=1\}$ without holomorphic points
satisfies $|g|<1$ and
\begin{eqnarray}
&&{\frac {G_{\bar z}}{|g|^{2}g_{\bar z}}}>0,\\
&&{\frac {|g|^{2}|G_{ z}|}{|G_{\bar z}|}}<1,\\
&&G_{z{\bar z}}+{\frac {{\bar g}_{ z}}{(|g|^{4}-1){\bar
g}}}G_{\bar z}-{\frac {|g|^{2}{\bar g}g_{\bar z}}{|g|^{4}-1}}G_{
z}=0,
\end{eqnarray}
put
\begin{eqnarray}
x_{1}&=& Re\left \{G-{\frac {1+|g|^{2}}{{\bar
g}g_{\bar z}}}G_{\bar z}\right\},\\
x_{2}&=& Im\left \{G-{\frac {1+|g|^{2}}{{\bar
g}g_{\bar z}}}G_{\bar z}\right\},\\
x_{3} &=& {\frac {1+|g|^{2}}{|g|^{2}g_{\bar z}}}G_{\bar z}.
\end{eqnarray}
Then $x=(x_{1},x_{2},x_{3}):M\rightarrow S_{1}^{3}$ is a space-like
surface with  de Sitter Gauss map $G$ and  antiholomorphic normal
Gauss map $g$ and Gauss curvature $K$ satisfying ${\sqrt
{1-K}}={\frac {1+|g|^{2}}{1-|g|^{2}}}.$ And the conformal structure
on $M$ is induced by the negative definite second fundamental form.
Conversely, any surface $x:M\rightarrow S_{1}^{3}$ with ${\sqrt
{1-K}}={\frac {1+|g|^{2}}{1-|g|^{2}}}(=-\eta_{3})$ can be given by
(5.16)(5.17)(5.18) and the de Sitter Gauss map $G$ and the normal
Gauss map $g$ must satisfy (5.13)(5.14)(5.15), where the conformal
structure on $M$ is induced by the negative definite second
fundamental form.}

\section{Graphs and examples}\setcounter{equation}{0}
In this section,we give the examples of  surfaces in $S_{1}^{3}$
with conformal normal Gauss maps within the  translational surfaces
and the Euclidean ruled surfaces.

In $H^{3}$, the graph $(u,v,f(u,v))$ with conformal normal Gauss map
satisfies the following fully nonlinear PDE (cf.[16][17])
\begin{equation}
f(f_{uu}f_{vv}-f_{uv}^{2})+[(1+f_{v}^{2})f_{uu}
-2f_{u}f_{v}f_{uv}+(1+f_{u}^{2})f_{vv}]=0.
\end{equation}

Take the  upper half-space model of  $S_{1}^{3}$. Consider the
space-like graph $(u,v,\\f(u,v))$ in $ S_{1}^{3}$ with  $
f_{u}^{2}+f_{v}^{2}<1.$ Its Gauss curvature  is given by $K=$
$$1-{\frac {f^{2}(f_{uu}f_{vv}-f_{uv}^{2})-f[(1-f_{v}^{2})f_{uu}
+2f_{u}f_{v}f_{uv}+(1-f_{u}^{2})f_{vv}]+(1-f_{u}^{2}-f_{v}^{2})}
{(1-f_{u}^{2}-f_{v}^{2})^{2}}}.$$ So $K=1-\eta_{3}^{2}$ is equvalent
to
\begin{equation}
f(f_{uu}f_{vv}-f_{uv}^{2})-[(1-f_{v}^{2})f_{uu}
+2f_{u}f_{v}f_{uv}+(1-f_{u}^{2})f_{vv}]=0,
\end{equation}
where $f_{u}^{2}+f_{v}^{2}<1.$ This is the fully nonlinear PDE which
the space-like graph in  $S_{1}^{3}$ with $K=1-\eta_{3}^{2}$ must
satisfy.\vspace{3mm}

{\sl Remark.} There exists a nice duality between the solutions of
minimal surface equation$$(1+f_{v}^{2})f_{uu}
-2f_{u}f_{v}f_{uv}+(1+f_{u}^{2})f_{vv}=0$$ in $R^{3}$ and the ones
of maximal surface equation$$(1-f_{v}^{2})f_{uu}
+2f_{u}f_{v}f_{uv}+(1-f_{u}^{2})f_{vv}=0$$ in Lorentz-Minkowski
3-space $L^{3}$(cf.[2]). Here, by the duality given by (5.2)(or
(5.3)), we know that if $f(u,v)$ is the solution of (6.1), then the
local graph of the space-like surface $(-ff_{u}-u,-ff_{v}-v,f{\sqrt
{1+f_{u}^{2}+f_{v}^{2}}})$  in  $ S_{1}^{3}$ satisfies (6.2).
Conversely, if $f(u,v)$ is the solution of (6.2) with
$f_{u}^{2}+f_{v}^{2}<1,$
 then the local graph of the  surface $(ff_{u}-u,ff_{v}-v,f{\sqrt
{1-f_{u}^{2}-f_{v}^{2}}})$  in $H^{3}$  satisfies (6.1).\vspace{3mm}

Next,as similar as done in section 6 of [16], we get the following
Theorem.\vspace{3mm}

\noindent{\bf THEOREM 6.1.} {\sl The nontrivial translational
space-like surfaces with the form $f(u,v)=\phi(u)+\psi(v)$ in
$S_{1}^{3}$ with conformal normal Gauss map are given, up to a
linear translation of variables, by
 \begin{equation}
 f(u,v)={\sqrt {a^{2}+u^{2}}}\pm {\sqrt {b^{2}+v^{2}}}
\end{equation}
with $f_{u}^{2}+f_{v}^{2}<1,$  where $a$ and $b$ are nonzero
constants.  The parameter form of these translational surfaces are
locally given by
\begin{equation}
x(u,v)=(a\sinh u,b\sinh v, a\cosh u+b\cosh v).
\end{equation}}

Considered as surfaces in 3-dimensional Minkowski space $L^{3}$, the
space-like ruled surfaces in $S_{1}^{3}$ can be represented as
$x(u,v)=\alpha(v)+u\beta(v):D\rightarrow S_{1}^{3},$ where
$D(\subset R^{2})$ is a parameter domain and $\alpha(v)$ and
$\beta(v)$ are two vector value functions into $L^{3}$ corresponding
to two curves in $L^{3}.$ When $\beta$ is locally nonconstant,
without loss of generality we can assume that either $\langle
\beta,\beta\rangle=1,$$\langle
\beta^{\prime},\beta^{\prime}\rangle=\pm 1,$ and $\langle
\alpha^{\prime},\beta^{\prime}\rangle=0$ or $\langle
\beta,\beta\rangle=1,$$\langle
\beta^{\prime},\beta^{\prime}\rangle=0,$ and $\langle
\alpha^{\prime},\beta\rangle=0,$ where $\langle \cdot,\cdot\rangle $
is the scalar product in  $L^{3}.$ As similar as done in Theorem 2
of [16], we have\vspace{3mm}

\noindent{\bf THEOREM 6.2.} {\sl Up to an isometric transformation
\begin{equation}
(x_{1},x_{2},x_{3})\rightarrow (x_{1}\cos \theta-x_{2}\sin
\theta+a,x_{1}\sin \theta+x_{2}\cos \theta+b, x_{3})
\end{equation}
 in $S_{1}^{3},$ every space-like ruled surface in
$S_{1}^{3}$ with conformal normal Gauss map is locally a part of one
of the following,\vspace{2mm}

(1) ordinary Euclidean space-like planes in $S_{1}^{3},$\vspace{1mm}

(2) $(u\cosh v,c\cdot\sinh v,u\sinh v),$ for a constant
$c\not=0,$\vspace{1mm}

(3) $(c_{2}\sinh v+u\cosh v,c_{1}\sinh v,c_{2}\cosh v+u\sinh v), $
for constant $c_{1}\not=0$ and $c_{2}\not=0.$} \vspace{3mm}

We should note that in the proof of Theorem 6.2, only when $\langle
\beta^{\prime},\beta^{\prime}\rangle =-1,$ we may get the nontrivial
cases (2) and (3).

Locally, the ruled surfaces (2) and (3) in Theorem 6.2  can be
represented as the  graph $(u,v,f(u,v))$ as follows,\vspace{3mm}

\noindent{\bf COROLLARY.} {\sl $f(u,v)=\pm {\frac
{c_{1}c_{2}+uv}{\sqrt {c_{1}^{2}+v^{2}}}}$ is  a solution of
equation (6.2), where $c_{1}\not=0$ and $c_{2}$ are
constants.}\vspace{3mm}

{\sl Remark.} In $H^{3},$ the translational surfaces
\begin{equation}
(a\cos u,b\cos v ,a\sin u+b\sin v)
\end{equation}
and the ruled surfaces
\begin{equation}
(u\cos v,c\cdot\sin v,u\sin v)
\end{equation}
and
\begin{equation}
(-c_{2}\sin v+u\cos v,c_{1}\cdot\sin v, c_{2}\cos v+u\sin v)
\end{equation}
with conformal normal Gauss map have been obtained ([16][17]), where
$a,b,c,c_{1}$ and $c_{2}$ are nonzero constants. Using (5.2) (or
(5.3)) and  Theorem 4.2, we may check that up to a isometric
transformation (6.5) in $S_{1}^{3}$ ($\theta=\pm {\frac {\pi}{2}}$),
(6.4) in Theorem 6.1 and (2) and (3) in Theorem 6.2 are,
respectively, the polar varieties of (6.6),(6.7) and (6.8) and vice
versa.\vspace{3mm}

{\sl Remark.} Every geodesic of $H^{3},$  corresponding respectively
to $u=0,$$u=\pi,$ $v=0$ and $v=\pi$ on surfaces (6.6) and to
$v={\frac {\pi}{2}}$ on surfaces (6.7) and to $v=\pm {\frac
{\pi}{2}}$ on surfaces (6.8) follow which $K=-1$ is mapped to a
simple point in $S_{0}$ by the generalized Gauss map.\vspace{3mm}

\section{Time-like surfaces in $S_{1}^{3}$ with conformal normal
Gauss map} In this section, we state the similar results as the
aboved for the time-like surfaces in $S_{1}^{3}$ without proofs.

Take the upper-half space model of $S_{1}^{3}$. Let $M$ be a
2-dimensional Lorentz surface and $x:M\rightarrow S_{1}^{3}$ be the
time-like immersiom with the local coordinates $u_{1},u_{2}.$ The
first and the second fundamental forms are given, respectively, by $
\textrm{I}=g_{ij}du_{i}du_{j}$ and $\textrm{II}=h_{ij}du_{i}du_{j}.$
The space-like unit normal vector is $N=x_{3}\eta_{1}{\frac
{\partial}{\partial x_{1}}}+ x_{3}\eta_{2}{\frac {\partial}{\partial
x_{2}}}+x_{3}\eta_{3}{\frac {\partial}{\partial x_{3}}},$  where
$\eta_{1}^{2}+\eta_{2}^{2}-\eta_{3}^{2}=1.$  Left-translating $N$ to
$T_{e}(R_{+}^{3}),$ we obtain
$$ \widetilde N:M\rightarrow S_{1}^{2}(1)\subset
T_{e}(R_{+}^{3}),$$
$$\widetilde N=L_{x^{-1}*}(N)=\eta_{1}{\frac {\partial}{\partial x_{1}}}(e)+
\eta_{2}{\frac {\partial}{\partial x_{2}}}(e)+ \eta_{3}{\frac
{\partial}{\partial x_{3}}}(e),$$ which is called the normal Gauss
map of time-like surface $x:M\rightarrow S_{1}^{3}$(cf.[1]).  Call
$\textrm{IV}=\langle {d\widetilde N},{d\widetilde N}\rangle$ the
fourth fundamental form of the time-like surface $x:M\rightarrow
S_{1}^{3}.$ We have
$\textrm{IV}=({\eta}_{3}^{2}g_{ij}-2\eta_{3}h_{ij}+g^{kl}h_{ik}h_
{jl})du_{i}du_{j}.$ Of course, we may also define the
high-dimensional version of the fourth fundamental form for the
time-like hypersurfaces in $S_{1}^{n+1}(1).$ \vspace{3mm}

\noindent{\bf THEOREM 7.1.} {\sl Let $M$ be a 2-dimensional Lorentz
surface and $x:M\rightarrow S_{1}^{3}$ be a time-like immersed
surface without umbilics. Then the normal Gauss map of $x(M)$ is
conformal if and only if the Gauss curvature
$K=1+\eta_{3}^{2}$.}\vspace{3mm}

In the Minkowski model of the de Sitter 3-space $S_{1}^{3},$ the
generalized Gauss map $N:M\rightarrow S_{1}^{3}$ of the time-like
surface $x:M\rightarrow S_{1}^{3}$ is a branched time-like immersion
with branch points where $K=1.$ \vspace{3mm}

\noindent{\bf THEOREM 7.2.} {\sl Let $M$ be a connected
2-dimensional Lorentz surface. Let $X:M\rightarrow S_{1}^{3}$ be a
time-like surface without umbilics and $K\not=1$. If the normal
Gauss map of $X:M\rightarrow S_{1}^{3}$ is conformal,
 then the normal Gauss map of its generalized Gauss map
 $N:M\rightarrow S_{1}^{3}$ is also
 conformal and vice versa.}\vspace{3mm}

The time-like graph $(u,v,f(u,v))$ in $S_{1}^{3}$ with conformal
normal Gauss map also satisfies the fully nonlinear PDE (6.2) with
$f_{u}^{2}+f_{v}^{2}>1.$\vspace{3mm}

\noindent{\bf THEOREM 7.3.} {\sl The nontrivial translational
time-like surfaces with the form $f(u,v)=\phi(u)+\psi(v)$ in
$S_{1}^{3}$ with conformal normal Gauss map are  given, up to a
linear transtation of variables, by\vspace{2mm}

(1) $f(u,v)={\sqrt {u^{2}+a^{2}}}\pm{\sqrt
{v^{2}+b^{2}}}$,\vspace{1mm}

(2) $f(u,v)={\sqrt {u^{2}-a^{2}}}\pm{\sqrt
{v^{2}-b^{2}}}$,\vspace{1mm}

(3) $f(u,v)={\sqrt {u^{2}+a^{2}}}\pm{\sqrt
{v^{2}-b^{2}}}$,\vspace{1mm}

(4) $f(u,v)={\sqrt {u^{2}-a^{2}}}-{\sqrt
{v^{2}+b^{2}}}$,\vspace{1mm}

(5) Flaherty time-like surface in $S_{1}^{3}$  (cf.[12]) $f(u,v)=\pm
u+\psi(v),$\vspace{2mm}\\
 where $a$ and $b$ are nonzero constants and
$\psi^{\prime}(v)\not=0$. }\vspace{3mm}

We may prove that the normal Gauss map of the time-like surfaces (2)
and (3) in Theorem 6.2 are also conformal. In addition, for the
time-like ruled surface  $x(u,v)=\alpha(v)+u\beta(v)$ in
$S_{1}^{3}$, we may also assume the remained four cases:\vspace{2mm}

(i) $\langle\beta,\beta\rangle
=-1,\langle\beta^{\prime},\beta^{\prime}\rangle=1,$ and $
\langle\alpha^{\prime},\beta^{\prime}\rangle=0,$\vspace{1mm}

(ii) $\beta$ is constant null vector,\vspace{1mm}

(iii) $\beta$ is constant and $\langle\beta,\beta\rangle =-1$,
$\langle\alpha^{\prime},\beta\rangle=0,$\vspace{1mm}

(iv) $\langle\beta,\beta\rangle =0,
\langle\beta^{\prime},\beta^{\prime}\rangle=1,$ and $
\langle\alpha^{\prime},\beta^{\prime}\rangle=0,$\vspace{2mm}\\
where $\langle\cdot,\cdot\rangle$ is the scalar product in $L^{3}$.
Hence, we have\vspace{3mm}

\noindent{\bf THEOREM 7.4.} {\sl Up to an isometric transformation
(6.5) in $S_{1}^{3},$ every time-like ruled surface in $S_{1}^{3}$
with conformal normal Gauss map is locally a part of one of the
following,\vspace{2mm}

(1) ordinary Euclidean time-like planes in $S_{1}^{3},$\vspace{1mm}

(2) ordinary Euclidean generalized cylinder
$x(u,v)=\alpha(v)+u\beta,$ where $\beta=(0,0,1)$ and $\alpha(v)$ is
arbitrary curve in $L^{3}$ with $\langle
\alpha^{\prime},\alpha^{\prime}\rangle>0,$\vspace{1mm}

(3) $(u\cosh v,c\cdot\sinh v,u\sinh v),$ for a constant
$c\not=0,$\vspace{1mm}

(4) $(c_{2}\sinh v+u\cosh v,c_{1}\sinh v,c_{2}\cosh v+u\sinh v), $
for constant $c_{1}\not=0$ and $c_{2}\not=0.$\vspace{1mm}

(5) $(u\sinh v,c\cdot\cosh v,u\cosh v),$ for a constant
$c\not=0,$\vspace{1mm}

(6) $(c_{2}\cosh v+u\sinh v,c_{1}\cosh v,c_{2}\sinh v+u\cosh v), $
for constant $c_{1}\not=0$ and $c_{2}\not=0,$\vspace{1mm}

(7) Flaherty's time-like surfaces in $S_{1}^{3}$ (cf[12]),
$x(u,v)=\alpha(v)+u\beta,$  where  $\beta=(1,0,1)$ and $\alpha(v)$
is arbitrary curve in $L^{3}$ with
$\langle\alpha^{\prime},\beta\rangle\not=0.$}\vspace{3mm}

We should note that in the proof of Theorem 7.4, only for case (i)
and (ii), we may get the surfaces (5)(6)(7) in Theorem 7.4. For case
(iv), we may assume $\beta(v)=(\rho(v)\cos \theta(v),\rho(v)\sin
\theta(v),\rho(v))$ with $\rho^{2}(\theta^{\prime})^{2}=1.$ Next, we
get a contradictory system of equations.\vspace{3mm}

{\sl Remark.} Up to a isometric transformation (6.5) in $S_{1}^{3}$
 $(\theta=\pm{\frac {\pi}{2}})$,  the time-like surfaces (3) and (4)
 in Theorem 7.4 are,  respectively, the polar varieties of the
 time-like surfaces (5) and (6) in Theorem 7.4 and vice versa.
 The similar result also holds for the time-like surfaces in
 Theorem 7.3. Generally, if $f(u,v)$ is the solution of (6.2) with
$f_{u}^{2}+f_{v}^{2}>1,$
 then the local graph of the  surface $(ff_{u}-u,ff_{v}-v,f{\sqrt
{f_{u}^{2}+f_{v}^{2}-1}})$  in $S_{1}^{3}$  also satisfies (6.2).
\vspace{3mm}

 Locally, the ruled surfaces (4) and (5) in Theorem 7.4 can be
represented as the graph $(u,v,f(u,v))$ as follows,\vspace{3mm}

\noindent{\bf COROLLARY.}  {\sl $f(u,v)=\pm {\frac
{c_{1}c_{2}-uv}{\sqrt {v^{2}-c_{1}^{2}}}}$ is  a solution of
equation (6.2), where $c_{1}\not=0$ and $c_{2}$ are
constants.}\vspace{3mm}

{\sl Remark.} When we do not assume that $f>0,$ (6.3) and
$f(u,v)=\pm {\frac {c_{1}c_{2}+uv}{\sqrt {c_{1}^{2}+v^{2}}}}$ and
$f(u,v)=\pm u+\psi(v),$ $ \psi^{\prime}(v)\not=0,$ are all
nontrivial entire solutions of the equation (6.2) defined on
$R^{2}.$ In addition, the cone $f(u,v)={\sqrt {u^{2}+v^{2}}}$ is
also the special solution of the equation (6.2), but its graph is
the light-like surface. By Omori-Yau's Maximum Principle[14][18],
there exist no entire solution  $f(u,v)$ of (6.2) satisfying
$f_{u}^{2}+f_{v}^{2}>1$ and $f>0$ on $R^{2}.$   Does there exist
nontrivial entire solutions of equation (6.2) defined on $R^{2}$
satisfying   $f_{u}^{2}+f_{v}^{2}<1$ and $f>0$?

\vspace{3mm}

\noindent  School of Mathematics and System Sciences,

 \noindent Shandong University,  Jinan 250100,

 \noindent P.R. China

\noindent E-mail: shishuguo@hotmail.com


\begin{thebibliography}{s16}
\bibitem{s1}Aiyama, R. and Akutagawa, K., Kenmotsu type
representation formula for space-like surfaces in the de Sitter
3-space, {\sl Tsukuba J. Math}, 24(2000), no.1, 189-196.
\bibitem{s2} Al{{\scriptsize {\'l}}}as, L.J. and Palmer, B., A
duality result between the minimal surface equation and the maximal
surface equation, {\sl Anais. Acad. Bras. Ci.}, 73(2001), no.2,
161-164.
\bibitem{s3}Bryant, R.L., Surfaces of mean curvature one in
hyperbolic space, {\sl Ast{\'e}risque}, 154-155(1987), 321-347.
\bibitem{s4}Epstein, C.L., The hyperbolic Gauss map and
quasiconformal reflections, {\sl J.Reine.Angew.Math.}, 372(1986),
96-135.
\bibitem{s5} G{\'{a}}lvez, J.A. and Mart{{\scriptsize {\'l}}}nez, A.,
The Gauss map and second fundamental form of surfaces in $R^{3}$,
{\sl Geom.Dedicata},  81(2000),181-192.
\bibitem{s6}G{\'{a}}lvez, J.A. and Mart{{\scriptsize {\'l}}}nez, A. and
Mil{\'{a}}n, F., Flat surfaces in the hyperbolic 3-space, {\sl
Math.Ann.}, 316(2000), no.3, 419-435.
\bibitem{s7}Klotz,T., Some uses of the second conformal structure on
strictly convex surfaces, {\sl Proc.Amer.Math.Soc.,} 14(1963),
793-799.
\bibitem{s8} Kokubu, M., Weierstrass representation for minimal
surfaces in hyperbolic space, {\sl T{\^{o}}hoku.Math.J.},49(1997),
367-377.
\bibitem{s9}Kokubu, M., Surfaces and fronts with harmonic-mean
curvature one in hyperbolic three-space, arXiv: math. DG/0504124.
\bibitem{s10} H. B. Lawson, Jr., Complete minimal surfaces in
$S^{3}$, {\sl Ann.Math.}, 92(1970), 335-374.
\bibitem{s11} Lee, S., Spacelike CMC 1 surfaces in de Sitter 3-space
$S_{1}^{3}(1)$: their construction and some examples, {\sl
Differential Geometry-Dynamical Systems}, 7(2005), 49-73.
\bibitem{s12}Milnor, T.K., A conformal analog of Bernstein's Theorem
for timelike surfaces in Minkowski 3-space, {\sl Contemporary
Mathematics}, Vol. 64, 1987, 123-132.
\bibitem{s13} Obata, M., The Gauss map of immersions of Riemannian
manifolds in spaces of constant curvature, {\sl J.Differential
Geom.}, 2(1968), 217-223.
\bibitem{s14} Omori, H., Isometric immersions of Riemannian
manifolds, {\sl J. Math. Soc. Japan.,} 19(1967), 205-214.
\bibitem{s15} Osserman, R., Minimal surfaces, Gauss maps, total
curvature, eigenvalue estimates, and stability, {\sl the Chern
Symposium 1979}, Springer, Berlin Heidelberg, New York 1980,
199-227.
\bibitem{s16} Shi, S.G., Weierstrass representation for surfaces
with prescribed normal Gauss map and Gauss curvature in $H^{3}$,
{\sl Chin.Ann.Math. B} , 25(2004), no.4, 567-586.
\bibitem{s17} Shi, S.G., New examples of surfaces in $H^{3}$ with
comformal normal Gauss map,  {\sl Anais. Acad. Bras. Ci.}, 78(2006),
no.1,7-16.
\bibitem{s18} Yau, S.T., Harmonic functions on complete Riemannian
manifolds, {\sl Comm. Pure. Appl. Math.}, 28(1975), 201-228.

\end{thebibliography}
\end{document}